\documentclass{amsart}

\usepackage{amssymb}
\usepackage[all]{xy}
\usepackage{hyperref}
\newtheorem{thm}{Theorem}%[section]

\newtheorem{lem}[thm]{Lemma}
\newtheorem{cor}[thm]{Corollary}

\newtheorem{prop}[thm]{Proposition}

\theoremstyle{definition}
\newtheorem{defn}[thm]{Definition}

\newtheorem{say}[thm]{}
\newtheorem{exmp}[thm]{Example}

\newtheorem{rem}[thm]{Remark}          

\newtheorem*{ack}{Acknowledgments}      % \renewcommand{\theack}{} 

\newtheorem{defn-thm}[thm]{Definition--Theorem}  %!!!!!!!!!!!!!!!!!!!!!!!!
\newtheorem{defn-lem}[thm]{Definition--Lemma}  %!!!!!!!!!!!!!!!!!!!!!!!!
\theoremstyle{remark}

\renewcommand{\c}[0]{{\mathbb C}}  

\renewcommand{\o}[0]{{\mathcal O}} 
\newcommand{\z}[0]{{\mathbb Z}}

\renewcommand{\r}[0]{{\mathbb R}} 

\renewcommand{\a}[0]{{\mathbb A}}

\newcommand{\p}[0]{{\mathbb P}}

\newcommand{\qtq}[1]{\quad\mbox{#1}\quad}
\newcommand{\spec}[0]{\operatorname{Spec}}
\newcommand{\pic}[0]{\operatorname{Pic}}

\newcommand{\red}[0]{\operatorname{red}}    
    
\newcommand{\im}[0]{\operatorname{im}}

\newcommand{\chr}[0]{\operatorname{char}}

\newcommand{\onto}[0]{\twoheadrightarrow}

\def\into{\DOTSB\lhook\joinrel\to}

\def\loccoh#1.#2.#3.#4.{H^{#1}_{#2}(#3,#4)}

\DeclareMathAlphabet{\mathchanc}{OT1}{pzc}%
                                {m}{it}

\usepackage[all]{xy}\xyoption{dvips}

\begin{document}
\bibliographystyle{amsalpha}
\hfill\today

  \title[Fundamental groups and path lifting]{Fundamental groups and path lifting\\ for algebraic varieties}
  \author{J\'anos Koll\'ar}

\begin{abstract} We study 3 basic questions about fundamental groups of algebraic varieties. For a morphism, is being   surjective on  $\pi_1$ preserved by base change?  What is the connection between openness in the Zariski  and in the Euclidean topologies? Which morphisms  have the  path lifting property?
\end{abstract}
  \maketitle

The aim of these notes is to study 3 questions involving maps between the 
fundemantal groups of algebraic varieties. 
\begin{itemize}
\item Let $X\to Y$ be a morphism of schemes that induces a surjection on the algebraic fundemantal groups. Does the same hold after a base change
$X\times_S Z\to Y\times_S Z$?
\item Let $X\to Y$ be a morphism between $\c$-schemes. When can we lift every continuous path in $Y(\c)$ to a path in $X(\c)$?
\item Let $X\to Y$ be a morphism between $\c$-schemes. What is the connection between openness in the Zariski  and in the Euclidean topologies?
\end{itemize}
An answer to the first question was used in the study of Pell surfaces
\cite{k-pell}. While this application involves only maps between algebraic curves, the curves in question are singular and non-proper, and it turns out to be not much harder to consider the general case.  This is treated in Section~\ref{sec.1}.

The proof uses some basic properties of open and universally open morphisms, some of which I did not find in the literature. These are worked out in 
Section~\ref{sec.2}.

Various forms of path lifting are studied in Section~\ref{sec.3}.
The answer is most complete for arc lifting (Definition~\ref{weak.l.p.say})
which is equivalent to  openness of the morphism in the Euclidean topology
and universal openness  in the Zariski topology; see Theorem~\ref{un.op.=.path.lift.loc.thm}.

While the main applications are to schemes of finite type, the discussions in
Sections~\ref{sec.1}--\ref{sec.2} are formulated for arbitrary    Noetherian  schemes.

\section{Maps between fundamental groups}\label{sec.1}

\begin{defn} \label{pi.surj.defn.1}
Let $X$ be a  connected scheme and $x\to X$ a geometric point. The {\it  fundamental group} of $X$ with base point $x$ is denoted by $\pi_1(X, x)$. Working with schemes, we use the 
algebraic fundamental group.

Let  $f:X\to Y$ be a  morphism of connected schemes.
Fix a base point $x\to X$ and its $f$-image  $y\to Y$. 
We get a natural group homomorphism
$f_*:\pi_1(X, x)\to \pi_1(Y, y)$.
\begin{enumerate}
\item We say that $f$ is {\it $\pi_1$-surjective} if
$f_*:\pi_1(X,x)\to \pi_1(Y,y)$
is surjective.
\end{enumerate}
Using the correspondence between quotients of the fundamental group and 
finite, \'etale covers we get the following equivalent form.
\begin{enumerate}\setcounter{enumi}{1}
\item  $f$ is $\pi_1$-surjective iff for every connected, finite, \'etale cover $Y'$, the fiber product   $X\times_YY'$ is also connected.
\end{enumerate}
\end{defn}
 The latter formulation shows that  the base point can be ignored in this definition.
More generally, choosing a different $x'\to y$, the image of 
$\pi_1(X, x')\to \pi_1(Y, y)$ is a conjugate of the 
image of 
$\pi_1(X, x)\to \pi_1(Y, y)$.

One of the main aims is to show that connectedness of the 
fiber product   $X\times_YY'$ also holds if $Y'\to Y$ is proper and universally open, see Theorem~\ref{pi1.onto.equiv.prop}.5.
We discuss open and universally open morphisms in Section~\ref{sec.2}.

\begin{thm} \label{pi1.onto.equiv.prop}
Let $X, S$ be   connected schemes and $g:X\to S$
  a proper and universally open morphism. The following are equivalent.
\begin{enumerate}
\item  $g$ is  $\pi_1$-surjective.
\item   $X\times_SX$ is  connected.
\item  $X\times_S\dots \times_SX$  ($n$ copies of $X$)  is  connected for some $n\geq 2$.
\item The number of  connected components of $X\times_S\dots \times_SX$ is bounded, independent of the number of factors.
\item For every    connected scheme $Y$ and 
   proper, universally open morphism $h:Y\to S$, the fiber product
$X\times_SY$ is  connected.
\item For every   connected scheme $Y$ and 
   proper, universally open morphism $h:Y\to S$, the second projection
$X\times_SY\to Y$ is  $\pi_1$-surjective.
\end{enumerate}
  \end{thm}

The following example shows that in
(\ref{pi1.onto.equiv.prop}.5--6)
we need  $Y\to S$ to be  universally open. It would not be   enough to assume only that $Y\to S$ is
 finite and open.

\begin{exmp}\label{pi1.onto.equiv.exmp}
 Let $C$ be a projective, nodal rational curve over an algebraically closed field 
with normalization $\pi:\p^1\cong \bar C\to C$ and $c\in C$ a smooth point.
Note that $\pi_1(C,c)\cong \z$.
Let $g_n:(c_n, C_n)\to (c, C)$ be its unique degree $n$, connected,  \'etale cover. For $n, m\geq 1$ 
set $(x,X_{n,m}):=(c_n, C_n)\amalg_{c_n=c_m}(c_m,C_m)$. 
The maps $g_n, g_m$ glue to 
$g_{n,m}:(x,X_{n,m})\to (c, C)$. 
Then
\begin{enumerate}
\item  $g_{n,m}$ is  $\pi_1$-surjective iff $(n, m)=1$,
\item $\bar C\to C$ is finite, open  and surjective, yet
\item $X_{n,m}\times_C\bar C$  is the disjoint union of
  $\bar C\amalg_{c=c}\bar C$ and of
 $n+m-2$ copies of $\bar C$. Thus it is disconnected iff $n+m>2$. 
\end{enumerate}
\end{exmp}

We start the proof of Theorem~\ref{pi1.onto.equiv.prop} with a series of remarks and lemmas that establish various special cases, and then use them to settle the general case.

\begin{rem}[Stein factorization]\label{stein.rem}
 Let $h:X\to X'$ be a proper morphism such that $h_*\o_X=\o_{X'}$. Then $X$ is connected if $X'$ is connected. 
Since $h_*$ commutes with any flat base change, (\ref{pi.surj.defn.1}.2) shows that $h_*: \pi_1(X, x)\to \pi_1(X', x')$ is an isomorphism.

Applying this to the Stein factorization $g:X\to X'\to S$ shows that  Theorem~\ref{pi1.onto.equiv.prop} holds for proper morphisms iff it holds for finite morphisms.
\end{rem}

\begin{lem}[\ref{pi1.onto.equiv.prop}.2$\Rightarrow$ \ref{pi1.onto.equiv.prop}.5] Let $X, Y, S$ be connected schemes and
  $g:X\to S, h:Y\to S$ finite, universally open morphisms. Assume that
  $X\times_SX$ is connected.  Then $X\times_SY$ is also  connected.
\end{lem}

Proof. 
Fix a geometric point $y\to Y$ and then choose $x_1, x_2\to X$   such that  $g(x_1)=g(x_2)=h(y)$.
Let $W\subset X\times_SX\times_SY$ be the connected component that contains
$(x_1, x_2, y)$. Since $h$ is finite and universally open, so is the
projection
$\pi_{12}:X\times_SX\times_SY\to X\times_SX$. Thus $\pi_{12}:W\to
X\times_SX$  is surjective by (\ref{basic.props.say}.3). In particular, there is a point $(x_1,
x_1, y')\in W$ for some $y'\to Y$.
Consider now the 2 projections  $\pi_i:W\to X\times_SY$ for $i=1,2$. Note
that
$$
(x_1, y), (x_1, y')\in \pi_1(W)\qtq{and}
(x_2, y), (x_1, y')\in \pi_2(W).
$$
Since $W$ is connected, this shows that $(x_1, y), (x_2, y)$ are in the same
connected component of $X\times_SY$.

Let $V_j\subset X\times_SY$ be the connected components.
Projection to $Y$ is universally open, so the projections $V_j\to Y$ are
surjective  by (\ref{basic.props.say}.3).
All preimages of $y\in Y$ are in the same connected component by the above argument, hence
$X\times_SY$ has only 1 connected component. \qed

\medskip
Using (\ref{pi.surj.defn.1}.2) this implies the following.

\begin{cor}[\ref{pi1.onto.equiv.prop}.2 $\Leftrightarrow$ \ref{pi1.onto.equiv.prop}.1]\label{pi1.onto.equiv.prop.cor.1}
  Let $X,  S$ be connected schemes and
  $g:X\to S$ a finite, universally open morphism. Assume that
  $X\times_SX$ is connected.  Then  $\pi_1(X)\to \pi_1(S)$  is surjective. \qed
  \end{cor}

\begin{lem} \label{etale.crit.lem}
Let $g:X\to S$ be a universally open morphism of finite type.  Assume that  the diagonal $\Delta_{X/S}$ is a connected component of $X\times_SX$. Then  $g$ uniquely factors as $g:X\to S'\to S$ where $X\to S'$ is a universal homeomorphism and $S'\to S$ is \'etale.
\end{lem}

Proof. Such a factorization is unique, so it is enough to construct it \'etale locally on $S$.
By Corollary~\ref{qf.to.f.by.et.lem}, after an \'etale base change we may assume that $g$ is of the form 
$g:\amalg_i (x_i,X_i)\to (s, S)$ where each restriction $g_i:(x_i,X_i)\to (s, S)$ is finite,  local and $k(x_i)/k(s)$ is purely inseparable.  
Thus  $X_i\times_SX_i$ is connected, since all of its irreducible components contain the 1-pointed scheme  $x_i\times_s x_i$. Thus
$\Delta_{X_i/S}=X_i\times_SX_i$,  so
  $k(x'_i)/k(g(x'_i))$ is also purely inseparable for every $x'_i\in X_i$. Hence $g_i$ is a universal homeomorphism by \cite[I.3.7--8]{ega71}.  Thus the factorization is  $g:\amalg_i (x_i,X_i)\to \amalg_i  (s, S)\to (s, S)$. \qed

\begin{lem} \label{etale.factor.lem}
Let  $X$ be a connected scheme and  $g:X\to S$  a finite, universally open morphism. Then  $g$ uniquely factors as $g:X\to S'\to S$
  where $S'\to S$ is finite, \'etale and $X\times_{S'}X$ is connected.
\end{lem}

Proof.  
Let $\Delta_{X/S}^{\rm conn}$ denote the  connected component of  $\Delta_{X/S}$ in $X\times_SX$.
It is a finite equivalence relation on $X$ and the geometric quotient $S_1:=X/\Delta_{X/S}^{\rm conn}$ exists by \cite[Lem.17]{MR2931872}. 
The natural map
$X\times_{S_1}X\to  X\times_SX$ is a universal homeomorphism onto  $\Delta_{X/S}^{\rm conn}$,
thus $X\times_{S_1}X$ is connected. 
We apply Lemma~\ref{etale.crit.lem} to $S_1\to S$ to get  $S_1\to S'\to S$
where  $S_1\to S'$ is a universal homeomorphism and $S'\to S$ is \'etale. 
\qed

\medskip

Combining Lemma~\ref{etale.factor.lem} and Corollary~\ref{pi1.onto.equiv.prop.cor.1}
 we get the finite case of following. For   proper morphisms  we also 
use  Remark~\ref{stein.rem}.

\begin{cor} \label{fin.ind.cor}
Let $g:X\to S$ be a proper, universally open morphism of  connected schemes. Let $x\to X$ be a geometric point and $s\to S$ its image. 
 Then
$\im[\pi_1(X, x)\to \pi_1(S, s)]$ has finite index in
$\pi_1(S, s)$. \qed
\end{cor}

The normalization of the nodal plane cubic shows that
the conclusion does not hold if $g$ is only assumed 
finite and open.

\begin{say}[Proof of Theorem~\ref{pi1.onto.equiv.prop}]\label{pi1.onto.equiv.prop.pf}
As we noted in Remark~\ref{stein.rem}, it is enough to prove the special case when $X\to S$ and $Y\to S$ are both finite. 

% Let $X\to S$ be a proper morphism and 
% $X\to X'\to S$ its Stein factorization. Then $X\to X'$ is $\pi_1$-surjective and $X$ is connected iff $X'$ is connected. Thus
% $X\to S$ is $\pi_1$-surjective iff $X'\to S$ is $\pi_1$-surjective.
% Applying this argument both to $X\to S$ and $Y\to S$, the proof of 
%  Theorem~\ref{pi1.onto.equiv.prop} is reduced to 

  We already proved that \ref{pi1.onto.equiv.prop}.1 $\Leftrightarrow$ \ref{pi1.onto.equiv.prop}.2 and that they imply 
\ref{pi1.onto.equiv.prop}.5.  Setting $X=Y$ shows that
\ref{pi1.onto.equiv.prop}.5 $\Rightarrow$ \ref{pi1.onto.equiv.prop}.2
and we get \ref{pi1.onto.equiv.prop}.3 by induction on the number of factors.

We use the shorthand $X^n_S:=X\times_S\dots \times_SX$  for the fiber product with $n$ factors.
The coordinate projection $\pi_{12}:X^n_S\to X\times_SX$
is surjective, thus \ref{pi1.onto.equiv.prop}.3 $\Rightarrow$ \ref{pi1.onto.equiv.prop}.2.

Since every
connected component of $X^n_S$  dominates $S$, if 
$X^n_S$  has at least $2$ connected components then
$X^{mn}_S$  has at least $2^m$ connected components. 
Thus (\ref{pi1.onto.equiv.prop}.3) $\Leftrightarrow$ (\ref{pi1.onto.equiv.prop}.4).

 If $g$ is not $\pi_1$-surjective then by Lemma~\ref{etale.factor.lem} it factors as
 $X\to S'\to S$ where $S'\to S$ is finite, \'etale and of degree $\geq 2$. 
Setting $Y:=S'$ shows  that (\ref{pi1.onto.equiv.prop}.6) $\Rightarrow$ (\ref{pi1.onto.equiv.prop}.1).

Conversely, assume (\ref{pi1.onto.equiv.prop}.1) and fix  $h:Y\to S$.
We already know that  $X\times_SY$ is connected. 
If $X\times_SY\to Y$ is  not $\pi_1$-surjective then by (\ref{pi.surj.defn.1}.2)
there is a nontrivial, finite \'etale cover  $Y'\to Y$ such that
$(X\times_SY)\times_YY'$ is not connected.
Applying (\ref{pi1.onto.equiv.prop}.5) to  $Y'\to S$ we get that
$X\times_SY'\cong X\times_SY\times_YY'$ is  connected, a contradiction. 
Thus (\ref{pi1.onto.equiv.prop}.1) $\Rightarrow$ (\ref{pi1.onto.equiv.prop}.6).
\qed
\end{say}

Next we consider 3 variants of $\pi_1$-surjectivity.

\begin{say}[Topological fundamental group] \label{pi.surj.defn.2.1}
Let $X$ be a connected $\c$-scheme of finite type and $x\in X$ a  point. 
We then have the {\it topological fundamental group}
$\pi^{\rm top}\bigl(X(\c),x\bigr)$. 
A morphism  $f:X\to Y$ between connected $\c$-schemes of finite type 
is   {\it $\pi^{\rm top}$-surjective} if 
$f_*: \pi^{\rm top}(X)\to \pi^{\rm top}(Y)$ is surjective.

Although the natural map $\pi^{\rm top}(X)\to \pi_1(X)$ can have infinite kernel, 
as a consequence of  
Lemma~\ref{fin.ind.cor.euc}
we see that if $f$ is proper and universally open  then the index of
$\im\bigl[\pi^{\rm top}(X)\to \pi^{\rm top}(Y)\bigr]$ in 
$\pi^{\rm top}(Y)$ equals the 
index of
$\im\bigl[\pi(X)\to \pi(Y)\bigr]$ in 
$\pi(Y)$. In particular,
$f$ is $\pi^{\rm top}$-surjective iff it is $\pi$-surjective.
Thus Theorem~\ref{pi1.onto.equiv.prop}  holds for the
topological fundamental group as well. 
\end{say}

\begin{say}[First homology group] \label{pi.surj.defn.2.2}
Let $X$ be a  connected scheme and $x\to X$ a geometric point. The {\it first homology group,}
denoted by $H_1(X)$, is defined as  the abelianization of $\pi_1(X,x)$; it is independent of the base point. Since we start with the algebraic fundamental group, $H_1(X)$
 is a $\hat\z$ module, where $\hat \z\cong \oplus_p \z_p$ is the profinite completion of $\z$. 
If $k=\c$ then   $H_1(X)$ is the  profinite completion of
$H_1\bigl(X(\c), \z\bigr)$.

We say that $f$ is  {\it $H_1$-surjective} if 
$f_*: H_1(X)\to H_1(Y)$
is surjective. Note that $\pi_1$-surjective
implies  $H_1$-surjective.

In \cite{k-pell} 
we needed to understand whether $H_1$-surjectivity 
is preserved by base change as in Theorem~\ref{pi1.onto.equiv.prop}.6.
The following example shows that it is not.

Let $X$ be a simply connected manifold (or variety over $\c$)  on which $A_n$ acts freely. Assume that $n\geq 6$ is odd. 
Let $A_{n-1}\subset A_n$ be a point stabilizer and  $C_n\subset A_n$  a subgroup generated by an $n$-cycle. We get a commutative diagram
$$
\begin{array}{ccc}
X & \stackrel{g'}{\longrightarrow} & X/C_n\\
\downarrow && \downarrow \\
X/A_{n-1} & \stackrel{g}{\longrightarrow} & X/A_n,
\end{array}
\eqno{(\ref{pi.surj.defn.2.2}.1)}
$$
which is a fiber product square.  Here
 $g$ is  $H_1$-surjective but  $g'$ is not.

This is the reason why, although in  \cite{k-pell}  the main interest is in 
$H_1$-surjectivity, we needed to 
understand the base change behaviour of $\pi_1$-surjectivity.
\end{say}

\begin{say}[Tame fundamental group] \label{pi.surj.defn.2.3}
Let $X$ be a $k$-scheme and $\chr k=p>0$. Let
$\pi_1^{(p)}(X, x) $ denote the largest prime to $p$ quotient of $\pi_1(X, x)$,
that is, the inverse limit of all quotients  $\pi_1(X, x)\onto H$ where  
$p\nmid |H|$.  
We say that $f$ is {\it $\pi_1$-surjective modulo $p$}
 if 
$f_*:\pi_1^{(p)}(X)\to \pi_1^{(p)}(Y)$ is surjective. 

The diagram  (\ref{pi.surj.defn.2.2}.1) also shows that $\pi_1$-surjectivity modulo $p$
is not preserved by base change. 
\end{say}

\section{Open and universally open maps}\label{sec.2}

We discuss properties of  universally open maps that were used in the proofs of
Theorems~\ref{pi1.onto.equiv.prop} and \ref{un.op.=.path.lift.loc.thm}.
We aim to treat these in their natural generality and also establish various results that are of independent interest.

\begin{defn}\label{op.uop.defn}
 A morphism $f:X\to S$ is {\it open at $x\in X$}
if for every open $x\in U\subset X$, its image $f(U)$ contains an open neighborhood of $f(x)$. 
A morphism $f:X\to S$ is {\it  open}  (resp.\ {\it open along $W\subset X$}) if $f$ is open at every $x\in X$  (resp.\ $x\in W$).

A morphism $f:X\to S$ is {\it universally open at $x\in X$}
(resp.\ along $W\subset X$)
if $f_T:X\times_ST\to T$ is open along $g_X^{-1}(x)$ (resp.\  $g_X^{-1}(W)$) for every $g:T\to S$ where $g_X:X\times_ST\to X$ is the first projection. 
We say that $f$ is {\it universally open} if it is  universally open at every $x\in X$. (Note that the Zariski topology on the product of 2 varieties is not the product topology, this is why open does not imply universally open.)

The following examples are good to keep in mind.
\begin{enumerate}
\item Let $g:X\to Y$ be a morphism of finite type $\c$-schemes. 
We see in Theorem~\ref{un.op.=.path.lift.loc.thm} that $g$ is universally open iff  $g(\c):X(\c)\to Y(\c)$ is open in the Euclidean topology.  Thus universal openness should be viewed as the more geometric notion.
\item If $g:X\to Y$ is open then  every irreducible component of $X$ dominates an  irreducible component of $Y$.  See Lemma~\ref{op.uop.defn.3.con} for a partial  converse statement. 
\item Let $g:C_1\to C_2$ be a quasi-finite morphism of purely 1-dimensional $k$-schemes.
If $C_2$ is irreducible then $g$ is open since every dense, constructible subset of $C_2$ is open. 
\item  Let $p:\bar X\to X$ be a finite, birational morphism. Then $p$ is universally open iff
it is a bijection on geometric points. 
\item As a consequence we see that if  $C$ is an irreducible  curve with nodes then the normalization $\bar C\to C$ is open but not universally open.
\item Let $p:\bar X\to X$ be a finite, birational morphism. Assume that every irreducible component of $X$ has dimension $\geq 2$. Then $p$ is open $\Leftrightarrow$ $p$ is  universally open $\Leftrightarrow$
$p$  is a bijection on geometric points. 
\end{enumerate}
Thus the difference between open and universally open appears mostly for 1-dimensional targets. Nonethless, many proofs involve localization and induction on the dimension, so the difference between the 2 notions can be significant.
\end{defn}

\begin{exmp}[Openness is not an open property]\label{op.not.op.exmp}
The following examples show that 
openness of a morphism at a point is not an open property. I assume this is why
this notion  is not defined in \cite[{Tag 004E}]{stacks-project}.
Nonetheless, I think that the notion is natural and has several useful properties.

(\ref{op.not.op.exmp}.1) Set $X:=(z=0)\cup(y=z-x=0)\subset \c^3$, $S:=(z=0)$
and $p:X\to S$ the coordinate projection. Then $p$ is universally open at
all points of the plane $(z=0)$ but not open at the points of
the punctured line $(y=z-x=0)\setminus\{(0,0,0)\}$. It truns out to be a general feature that openness along a fiber says a lot about the maximal dimensional irreducible components of $X$ but very little about the lower 
dimensional ones; see Theorem~\ref{uo.pured.sub.thm}.

(\ref{op.not.op.exmp}.2) 
 Let $Y$ be the pinch point
$(x^2=y^2z)\subset \c^3$ and $X\cong \c^2_{uv}$ its normalization.
Here $p:X\to Y$ is given by  $(u,v)\mapsto  (uv, u, v^2)$. Note that
$g$ is universally open at $(0,0)$ (see (\ref{vol.op.crit}.3) for a more general claim).
However it is not open at any other point of the $v$-axis.

Furthermore, although all fibers of $p:X\to Y$ have dimension 0,
it does not have pure relative dimension 0. Indeed, the fiber product
$X\times_YX$ has 2 irreducible components; one is $X\cong \c^2$ and the other
is isomorphic to $\c^1$ lying over the $v$-axis.

See Example~\ref{op.not.op.exmp.2} for other properties  of this surface.

\end{exmp}

\begin{say}[Basic properties]\label{basic.props.say}
The following are some obvious properties.
\begin{enumerate}
\item  Let  $f:X\to S$ be a morphism and $Z\subset X$ a locally closed subscheme.
If $f|_Z$ is open at $x\in Z$ then $f$ is open at $x$.
\item Let $g:Y\to X$ and $f:X\to S$ be morphisms,  $y\in Y$. 
\begin{enumerate}
\item
If $g$ is open at $y$ and $f$ is open at $g(y)$ then  $f\circ g$ is open at $y$.
\item If $f\circ g$ is open at $y$ then $f$ is open at $g(y)$. 
\end{enumerate}
\item  %(\ref{univ.op.exmp}.1)
Let $g:X\to Y$ be proper and open. Let $X_i\subset X$ be a connected component. Then $g(X_i)$ is closed and open in $Y$, hence a connected component.
\item If  $f:X\to S$ is  open at $x\in X$ then  $f_T:f^{-1}(T)\to T$ is also open at $x$ for every subvariety
$f(x)\in T\subset S$.
\item Let  $f:X\to S$ be morphism and  $S_i\subset S$ the irreducible components. Set  $X_i:=f^{-1}(S_i)$. By restriction we get $f_i:X_i\to S_i$.
Then $f$ is open (resp.\ universally open) at $x$ iff $f_i$ is open (resp.\ universally open) at $x$ whenever $x\in X_i$. 
\item A flat morphism of finite presentation is universally open
\cite[Tag 01UA]{stacks-project}. This is probably easiest to see using (\ref{vol.op.crit}.3).
\item \cite[Tag 04R3]{stacks-project}
Being  universally open is \'etale local on source and target. That is, if we have a commutative diagram 
 $$
\begin{array}{ccc}
(x',X') & \stackrel{h'}{\longrightarrow} & (x,X)\\
g'\downarrow\hphantom{g'} && \hphantom{g}\downarrow g \\
(s',S') & \stackrel{h}{\longrightarrow} & (s,S),
\end{array}
%\eqno{(\ref{qf.to.f.by.et.lem}.1)}
$$
where $h,h'$ are \'etale, then $g$ is universally open at $x$ iff
 $g'$ is universally open at $x'$.
\item  Let $X\to S$ be a proper morphism with Stein factorization
$X\to X'\to S$. By (2.b), if  $X\to S$ is open then so is  $X'\to S$.
The following example shows that $X\to X'$ need not be open.
 Take 2 copies of  $\a^1\times \p^1$ and  and glue them along the points $(0_i, p)\in \a^1_i\times \p^1$ for some $p\in \p^1$ to get $X$. 
The Stein factorization of the first coordinate projection is
$$
\pi:(\a^1_1\times \p^1)\amalg_{(0_1,p)=(0_2,p)}(\a^1_2\times \p^1)
\stackrel{\rho}{\to} 
\a^1_1\amalg_{0_1=0_2}\a^1_2\stackrel{\sigma}{\to} \a^1.
$$
Note that  $\pi$  and $\sigma$ are universally open but $\rho$ is not open.
\end{enumerate}
\end{say}

\begin{say}[Valuative and base change criteria]\label{vol.op.crit}
 The simplest non-open morphism is the
embedding of a closed point into an irreducible curve  $\{p\}\into C$. 
It turns out that this example is quite typical. 
 If there is an irreducible, positive dimensional subvariety
$s\in T\subset S$ such that $f^{-1}(T)=f^{-1}(s)$ then, by (\ref{basic.props.say}.4),
$f^{-1}(T)\to \{s\}\into T$ shows that $f$ is not open at $x$,
giving a necessary openness criterion.

We claim that if $f:X\to S$ is 
of finite type, then the criterion is also sufficient, after a  small change.

\medskip
(\ref{vol.op.crit}.1) Let $f:X\to S$ be a  morphism 
of finite type,  $x\in X$ a  point and $s:=f(x)$. Then $f$ is not open at $x$ iff
there is an open subset  $x\in U\subset X$ and an   irreducible   subscheme
$s\in C\subset S$ such that 
$s$ has codimension 1 in $C$ and the generic point of $C$ is not contained in 
$f(U)$.
In particular,    
$$
f|_{ U\cap f^{-1}(C)}:  U\cap f^{-1}(C)\to \bar s\to C
\eqno{(\ref{vol.op.crit}.1.a)}
$$
shows that $f$ is not open at $x$.
\medskip

Proof. Choose $U$ such that
$f(U)$ does not contain any open neighborhood of $s$. Since $f(U)$ is constructible, $S\setminus f(U)$ is constructible and its closure contains $s$. There is thus an irreducible component $W\subset S\setminus f(U)$
whose closure contains $s$.  Note that $f(U)\cap \bar W$ is a nowhere dense
constructible subset. Take any irreducible  subscheme
$\{s\}\in C\subset S$ that is not contained in the closure of
$f(U)\cap \bar W$.   \qed 
\medskip

Since $\o_{s,C}$ is dominated by a valuation ring, we can 
 restate (\ref{vol.op.crit}.1) in the following variant forms.

\medskip
(\ref{vol.op.crit}.2)   A finite type morphism $f:X\to S$  is open at $x\in X$
if $f_T:X\times_ST\to T$ is open along $g_X^{-1}(x)$ for every 1-dimensional, irreducible subscheme  $s\in T\subset \spec\o_{s,S}$ where $s:=f(x)$.  \qed

\medskip
(\ref{vol.op.crit}.3)  \cite[{Tag 01TZ}]{stacks-project}. A finite type morphism $f:X\to S$  is universally open at $x\in X$
if $f_T:X\times_ST\to T$ is open along $g_X^{-1}(x)$ 
 for every $g:T\to S$ where  $T$ is the spectrum of a valuation ring and $g$ maps its closed point to $s:=f(x)$.  \qed

\medskip
(\ref{vol.op.crit}.4)
  Let $g:X\to S$ be a morphism. Assume that there is a  closed subscheme  $x\in Z\subset X$ such that $h:=g|_Z:Z\to S$ is finite, dominant,   $h^{-1}(s)=\{x\}$  and $k(x)/k(s)$ is purely inseparable.  Then $g$ is universally open at $x$.

\medskip
Proof. Let  $x\in U\subset X$ be open. Then $Z\setminus U$ is closed, hence so is
$h(Z\setminus U)$, which  does not contain $s$. So
$$
s\in S\setminus h(Z\setminus U)\subset h(Z\cap U)\subset g(U)
$$
shows that $g$ is open at $x$. The assumptions are preserved by base change, so
$g$ is universally open at $x$. \qed

\medskip
(\ref{vol.op.crit}.5) Let  $g:X\to (s,S)$ be a morphism  of finite type.
Assume that $X,S$ are integral, $g$ is open along $X_s$. Then 
$\dim X=\dim X_s+\dim S$. 
\medskip

Proof. We need to prove that the generic fiber also has dimension $X_s$.
This is clear after base change to $g:T\to S$ as in (\ref{vol.op.crit}.3),
where $g$ maps the closed point to $s$ and the generic point of $T$ to the generic point of $S$. \qed

\medskip
(\ref{vol.op.crit}.6)
We prove in Theorem~\ref{uop.et.test.thm} that 
a finite type morphism $f:X\to (s,S)$  is universally open along $X_s$
iff $f_{S'}:X\times_SS'\to (s',S')$ is open along $(X\times_SS')_{s'}$ 
 for every \'etale $g:(s',S')\to (s,S)$.
\end{say}

\begin{say}[Openness and pure dimensional morphisms]\label{pd.crit.say}
Let  $g:X\to S$ be a morphism  of finite type.

\medskip
(\ref{pd.crit.say}.1) Set $ X^{(n)}:=\{x\in X: \dim_x X_{g(x)}=n\}$ and
$ X^{(\leq n)}:=\{x\in X: \dim_x X_{g(x)}\leq n\}$. 

By the upper semicontinuity of the fiber dimension \cite[Tag 02FZ]{stacks-project}, $X^{(\leq n)}$ is open in $X$ and
$ X^{(n)}=X^{(\leq n)}\setminus X^{(\leq n-1)}$ is  closed in $X^{(\leq n)} $ and locally closed in $X$.

\medskip
(\ref{pd.crit.say}.2) Let  $g:X\to (s,S)$ be a morphism  of finite type.
 Then
$g$ is open (resp.\ universally open) along $X_s$  iff  $g^{(n)}: X^{(n)}\to S$ is 
open (resp.\ universally open) along  $X^{(n)}_s$ for every $n$.

\medskip
Proof.   As we noted,  $X^{(\leq n)}$ is open in $X$, thus we may as well assume that
all fibers have dimension $\leq n$.
The formation of  $X^{(n)}$ commutes with base change, and over a 1-dimensional base the claims are clear. \qed
\medskip

Note that the punctual version does not hold. As an example, set $S:=(xy=0)\subset \a^2$ and $X=(x=0)\cup(y=z=0)\subset \a^3$. The coordinate projection is open at the origin yet $X^{(1)}\to S$ is not  open at the origin.
A much worse example is given in Example~\ref{op.n.uop.unibr.exmp}.

\medskip
(\ref{pd.crit.say}.3) Let  $g:X\to S$ be a morphism  of finite type whose fibers have pure dimension $n$.
Let $x\in Z\subset X$ be a relative complete intersection. 
Then $g$ is  open (resp.\ universally open) at $x$ iff $g|_Z$ is open (resp.\ universally open) at $x$.

\medskip
Proof. Assume that  $g|_Z$ is open at $x$ and let  $x\in U\subset X$ be open.
Then $Z\cap U\subset Z$ is open and $g(x)\in g(Z\cap U)\subset g(U)$ 
shows that $g$ is  open at $x$.

Conversely, assume that  $g|_Z$ is not open at $x$.
By (\ref{vol.op.crit}.1--2) it is enough to check the claim when
$S$ is a spectrum of a local ring of dimension 1.  Thus, after 
replacing $X$ with  an open neighborhood of $x$ we may asume that
$Z_s=Z$; this follows from (\ref{vol.op.crit}.1).
Since $\dim_xZ=\dim_xX-r$ and $\dim_xZ_s=\dim_xX_s-r$, we conclude that
$\dim_xX=\dim_xX_s$. Thus $X_s$ is an irreducible component of $X$, hence
$g$ is not open. \qed

\medskip
(\ref{pd.crit.say}.4)  Let  $g:X\to S$ be a morphism  of finite type.
If $g$ has {\it pure relative dimension} $n$ (that is,  $X\times_ST$ has pure dimension  $\dim T+n$ for every local morphism $T\to S$ from 
 an  integral scheme  to $S$) then $g$ is also universally open. 
This follows from (\ref{vol.op.crit}.3). 
Note, however, that the local version of this is not true; see (\ref{op.not.op.exmp}.2). 
\end{say}

Next we  show that it is enough to use 
\'etale base changes in the definition of universal openness.

\begin{thm}\label{uop.et.test.thm}
 Let  $g:X\to (s,S)$ be a morphism  of finite type.
Then $g$ is universally  open  along $X_s$ iff   $g'$ is open along $X'_{s'}$ for  every \'etale base change diagram as in (\ref{basic.props.say}.7). 
\end{thm}

Proof.  If $g$ is universally  open  at $x$ then it  is open after every  base change. 

Conversely, pick $x\in X_s$ and assume that $g$ is open at $x$ after every \'etale base change. Set $n=\dim_x X_s$. 
By  (\ref{pd.crit.say}.2)
 $g^{(n)}: X^{(n)}\to S$ is also  open at $x$ after every \'etale base change
and by (\ref{pd.crit.say}.3) the same holds for every relative complete intersection  $x\in Z\subset X^{(n)}$ of codimension $n$. Then  $g|_Z$ is quasi-finite.
Take an \'etale base change  $(s',S')\to (s, S)$ as in 
Proposition~\ref{qf.to.f.by.et.thm} to get
$$
Z\times_SS'=W\coprod \amalg_i  (z'_i, Z'_i).
$$
By construction every $ (z'_i, Z'_i)\to (s',S')$ is finite and
$k(z'_i)/k(s')$ is purely inseparable.   
Furthermore, $ (z'_i, Z'_i)\to (s',S')$ is open by assumption, hence dominant.
Thus 
$ (z'_i, Z'_i)\to (s',S')$ is  universally  open by (\ref{vol.op.crit}.4)
and so is $g|_Z:Z\to S$. Thus  $g$ is universally  open by (\ref{pd.crit.say}.3) and (\ref{pd.crit.say}.2). \qed

\medskip

We used some results  on \'etale localization of quasi-finite morphisms;
see \cite[{Tag 04HF}]{stacks-project} for proofs and further generalizations.

\begin{prop}\label{qf.to.f.by.et.thm}  
Let $g:X\to S$ be a quasi-finite morphism and $s\in S$.
Then there is an  \'etale morphism  $(s',S')\to (s, S)$ such that
$$
X\times_SS'=W\coprod \amalg_i\  (x'_i, X'_i),
$$
where, $W$ does not have any points lying over $s'$ and, for every $i$,  the morphism    $g'_i:(x'_i,X'_i)\to  (s',S')$
is finite,  $(g'_i)^{-1}(s')=\{x'_i\}$ and  $k(x'_i)/k(s')$ is purely inseparable.  \qed
\end{prop}

\begin{cor}\label{qf.to.f.by.et.lem} %\cite[{Tag 04HF}]{stacks-project}
Let $g:(x,X)\to (s,S)$ be a quasi-finite morphism.
Then  there is a commutative diagram
 $$
\begin{array}{ccc}
(x',X') & \stackrel{h'}{\longrightarrow} & (x,X)\\
g'\downarrow\hphantom{g'} && \hphantom{g}\downarrow g \\
(s',S') & \stackrel{h}{\longrightarrow} & (s,S),
\end{array}
\eqno{(\ref{qf.to.f.by.et.lem}.1)}
$$
where $h, h'$ are \'etale, $g'$ finite, 
${g'}^{-1}(s')=\{x'\}$ and  $k(x')/k(s')$ is purely inseparable. \qed
\end{cor}

We refer to (\ref{qf.to.f.by.et.lem}.1) as an {\it \'etale base change diagram} of $g$. 

\begin{cor}\label{qf.to.f.by.et.lem.2} \cite[{Tag 02LO}]{stacks-project}
Let $g:(x,X)\to (s,S)$ be a finite morphism.
Then there is an  \'etale morphism  $(s',S')\to (s, S)$ such that
$$
X\times_SS'=\amalg_i\  (x'_i, X'_i),
$$
where, for every $i$,  the morphism    $g'_i:(x'_i,X'_i)\to  (s',S')$
is finite,  $(g'_i)^{-1}(s')=\{x'_i\}$ and  $k(x'_i)/k(s')$ is purely inseparable.  \qed
\end{cor}

\subsection*{Geometrically  unibranch schemes}{\ }

\begin{defn}\label{unibr.defn} A local scheme $(x, X)$ is called {\it geometrically  unibranch} if for every \'etale morphism
$(x',X')\to (x, X)$, the  local scheme $(x', X')$ is irreducible.
See \cite[{Tag 06DT}]{stacks-project} for other definitions and basic properties. 
\end{defn}

Openness is very well behaved for geometrically  unibranch targets; cf.\ \cite[Tag 0F32]{stacks-project}

\begin{prop}\label{new.open.lem.2}
  Let $g:(x, X)\to (s, S)$ be a finite type morphism. Assume that $X$ is irreducible and 
 $S$ is geometrically  unibranch at $s$. 
The following are equivalent.
\begin{enumerate}
\item $\dim X=\dim X_s+\dim S$.
\item  $g$ is open along $X_s$.
\item   $g$ is universally open along $X_s$.
\end{enumerate}
\end{prop}

Proof.  Note that (3) implies (2)  by definition, and
(2) implies (1) by (\ref{vol.op.crit}.5).
 It remains to show that (1) $\Rightarrow$ (3).

Note that $\dim_x X\leq\dim_x X_s+\dim_s S$ for every $x\in X_s$, thus
$X_s$ is pure dimensional. 

Let $x\in Z\subset X$ be a relative complete intersection  of codimension $\dim X_s$. 
By (\ref{pd.crit.say}.3)   it is enough to show that
$g|_Z$ is  universally open at $x$. 
Thus  we may assume that $g$  is quasi-finite. Apply 
Proposition~\ref{qf.to.f.by.et.thm} to get an  \'etale morphism  $(s',S')\to (s, S)$ and a decomposition
$$
X\times_SS'=W\coprod \amalg_i\  (x'_i, X'_i),
$$
where the   $k(x'_i)\supset k(s')$ are purely inseparable and 
the projections  $g'_i:X'_i\to S'$ are finite. The $g'_i$ are  also dominant by assumption, hence  universally open at $x'_i$ by (\ref{vol.op.crit}.4). Thus $X\times_SS'$ is  universally open along  $X'_{s'}$, and so  $X\to S$ is  universally  open along $X_s$. \qed
\medskip

\begin{exmp} \label{op.n.uop.unibr.exmp}
%\begin{cor}\label{new.open.lem.1}
It seems natural to hope that the equivalence 
(\ref{new.open.lem.2}.2) $\Leftrightarrow$ (\ref{new.open.lem.2}.3)
holds pointwise. This is, however, not the case. To see this, we construct below a projective  morphism of surfaces  $g:(x, X)\to (s, S)$ such that 
 $S$ is normal and  $g$ is  open at $x$, yet $g$  is  not universally open at $x$.

 Let $(s,S)$ be a normal surface singularity with a (non-minimal) resolution  $\tau:Y\to S$ and exceptional curves $E_1,\dots, E_n$. 
Assume that
\begin{enumerate}
\item for every $i=1,\dots, n-1$ there is a morphism  $Y\to (x_i, X_i)$ that contracts only the curve $E_i$,
\item $X_i$ is normal and $\pi_i:X_i\to S$ is projective,
\item if an algebraic curve $C\subset Y$ is disjoint from 
 $E_1\cup\cdots\cup E_{n-1}$ then it is also disjoint from $E_n$.
\end{enumerate}
Let $(x,X)$ be obtained from the surfaces  $X_1,\dots, X_{n-1}$ by identifying the points $x_1,\dots, x_{n-1}$. The morphisms $\pi_i$ glue to a
projective morphism  $\pi:(x,X)\to (s,S)$.

\medskip
{\it Claim \ref{op.n.uop.unibr.exmp}.4.} $\pi$ is open but not universally open at $x$.
\medskip

Proof.  Let $B\subset Y$ be an algebraic curve that meets $E_n$ transversally at a general point $b\in E_n$. Then $\spec\o_{b,B}\to S$ shows that 
  $\pi$ is  not universally open at $x$.

Assume to the contrary that  $\pi$ is not open. Then, by (\ref{vol.op.crit}.1), there is a curve $0\in C_S\subset S$ such that the closure $C_X$ of the preimage
of $C_S\setminus\{s\}$ does not pass through $x$. Note that 
$C_X$ is the union of the birational transforms $(\pi_i)^{-1}_*C_S$
and these in turn are the images of $C:=\tau^{-1}_*C_S\subset Y$. 
Thus $C$ is disjoint from 
 $E_1\cup\cdots\cup E_{n-1}$. By (3) then $C$  is also disjoint from $E_n$,
which contradicts  $0\in C_X$. \qed
\medskip

In order to construct such an $(s,S)$, we start with  $\p^2$ and 3 general lines $L_1,L_2, L_3$. Let $B\subset \p^2$ be a  general quartic curve.
We obtain $Y'$ by blowing up the 12 intersection points of $B$ with the lines.
The  birational transforms of the $L_i$ become $E'_i\subset Y'$.
Contracting them we get $(s,S)$. 
Pick next a very general point $p\in E'_3$ and let $Y=B_pY'$ denote the blow-up. The  birational transforms of the $E'_i$  give $E_1, E_2, E_3$ and
$E_4$ is the exceptional curve of $Y\to Y'$. 

Each $E_i\subset Y$ is a rational curve with negative self-intersection, so it can be contracted projectively.
(This is essentially due to Castelnuovo; the proof in \cite[V.5.7]{hartsh} is easy to modify.) It remains to check assumption (3).
Slightly stronger, we claim that an algebraic curve $C'\subset Y'$ can not intersect   $E':=E'_1\cup E'_2\cup E'_3$ only at $p$. To see this note that
 $\pic(Y')$ is finitely generated  (in fact, isomorphic to $\z^{13}$), hence
the image of the restriction map 
$$\pic(Y')\to \pic(E')\cong \z+\c^*
$$
is finitely generated. Thus, for a very general point $p\in E'$, the intersection of 
$\z[p]\subset \pic(E')$ with the image is the trivial  element $[0]\in \pic(E')$. 

As a concrete example, we can take $(s,S)$ to be 
the projectivisation of the affine surface
$$
(s, S^0):=\bigl(xyz+x^4+y^4+z^4=0\bigr)\subset \a^3. \qed
$$
\end{exmp}

If $X$ is reducible and $g:X\to (s,S)$ is universally open along $X_s$
then, as shown by Example~\ref{prop.pl.exmps}.2,  we can say very little about the lower dimensional irreducible components of $X$. The next result  shows that the maximal dimensional irreducible components behave better.

\begin{thm} \label{uo.pured.sub.thm}
Let  $g:X\to (s, S)$  a finite type morphism that is  universally open along $X_s$. Set $n:=\dim X_s$ and let $X^{\rm max}\subset X^{(n)}$ be the union of  all those  irreducible components of $X^{(n)}$ that dominate some irreducible component of $S$ and have 
nonempty intersection with $X_s$.

Then $X^{\rm max}\to S$ is universally open  along $X_s\cap X^{\rm max}$.
\end{thm}

Proof. The conclusions can be checked after an \'etale base change, hence, using Lemma~\ref{unibr.after.etale},
we may assume that every  irreducible component of $(s,S)$ is unibranch.
Then, by (\ref{basic.props.say}.7), we may also assume that $(s,S)$ is irreducible, hence unibranch.

Let  $X_i\subset X^{\rm max}$ be an irreducible component. By assumption
$X_i$ dominates $S$ and the generic fiber of  $X_i\to S$ has dimension $n$ since
$X_i\subset X^{(n)}$. Thus $\dim X_i=n+\dim S$ and so   $X_i\to S$ is
 universally open along $X_s\cap X_i$ by Proposition~\ref{new.open.lem.2}.  \qed

\begin{lem}\cite[\href{https://stacks.math.columbia.edu/tag/0CB4}{Tag 0CB4}]{stacks-project}\label{unibr.after.etale}
 Let $(s, S)$ be a  local scheme.
Then there is an  \'etale morphism
$(s',S')\to (s, S)$ such that every irreducible component
$(s', S'_i)\subset (s',S')$ is geometrically  unibranch. \qed
\end{lem}

\medskip

The following is  a partial converse to (\ref{op.uop.defn}.3). 

\begin{lem}\label{op.uop.defn.3.con} Let $S$ be a connected scheme and $g:X\to S$  a dominant morphism of finite type whose fibers have pure dimension $n$. The following are equivalent.
\begin{enumerate}
\item $g:X\to S$ is universally open,
\item  $g\times g:X\times_SX\to S$ is open,
\item  Every irreducible component of $X$ dominates an
 irreducible component of $S$ and  every irreducible component of $X\times_SX$ dominates an irreducible component of $S$.
\end{enumerate}
\end{lem}

Proof. If  $g:X\to S$ is universaly open   then so is
$X\times_SX\to X$, hence also their composite $X\times_SX\to S$.
Thus  (1) $\Rightarrow$ (2) and  (2) $\Rightarrow$ (3) is clear. 

It remains to show that if (3) holds then $g$ is universally open along $X_s$ for every $s\in S$.
 Using Lemma~\ref{unibr.after.etale}, after an \'etale base change we may assume that every irreducible component
$(s,S_i)\subset (s,S)$ is geometrically  unibranch.

Set $X_i:=g^{-1}(S_i)$. By (\ref{basic.props.say}.7)  it is enough to show that  each $X_i\to S_i$  is universaly open along $(X_i)_s$.
If this does not hold, then, by
Proposition~\ref{new.open.lem.2}, there is an irreducible component
$Z_i\subset X_i$ that does not dominate $S_i$. 
By assumption there is an irreducible component  $X_j\subset X$
that contains $Z_i$ and this $X_j$ dominates some irreducible component  $S_j\subset S$. By assumption $S_i\neq S_j$. 

Let $X^*_i\subset X$ be an irreducible component that does  dominate $S_i$. 
Then $X_j\times_S X^*_i$ is an union of irreducible components of $X\times_SX$ that lie over  $S_i\cap S_j$, hence they do not dominate any
irreducible component of $S$.  This is impossible by (3), \qed

\begin{exmp} This example shows that the equidimensionality assumption in the previous 
Lemma~\ref{op.uop.defn.3.con}  is necessary.
Set
$$
X:=(x_0y_0+x_1y_1+x_2y_2=0)\subset \a^3_{\mathbf x}\times  \p^2_{\mathbf y}
$$
with projection $\pi:X\to S=\a^3_{\mathbf x}$. The central fiber has dimension 2, so $\pi$ is not  open.  
The fiber product  
$$
X\times_SX\subset \a^3_{\mathbf x}\times  \p^2_{\mathbf y}\times  \p^2_{\mathbf z}
$$ is given by 2 equations, hence its irreducible components have dimension
$\geq 7-2=5$. The   central fiber of $X\times_SX\to S$ has dimension 4,
so it can not be an irreducible component. Thus $X\times_SX$ is irreducible.

Note also that if $X\to S$ is not pure dimensional then the fiber product
$X\times_S\cdots\times_S X$ of more than $\dim S$ copies of $X$ always has a non-dominant 
irreducible component.
\end{exmp}

\section{Path lifting in the Euclidean topology}\label{sec.3}

Let  $g:X\to Y$ be a  morphism of $\c$-schemes of finite type. 
In this section we compare scheme-theoretic properties of $g$
with  properties of $g(\c):X(\c)\to Y(\c)$ in the Euclidean topology. 
 
It is easy to see that universal openness in the Zariski topology (Definition~\ref{op.uop.defn}) is equivalent to
openness in the Euclidean topology; see Lemma~\ref{vol.op.crit.7}.
Next we study 3 versions of the path lifting property. 
We see in Theorems~\ref{un.op.=.path.lift.loc.thm} and~\ref{un.op.=.path.lift.thm} that 2 of them have very satisfactory scheme-theretic descriptions.

\begin{defn}[Arc and path lifting]\label{weak.l.p.say}
  Let  $h:M\to N$  be a continuous map of topological spaces
and $\gamma:[0,1]\to N$ a path. A {\it lift} of
$\gamma$ with starting point $m\in h^{-1}\bigl(\gamma(0)\bigr)$ is a 
continuous map  $\gamma':[0,1]\to M$ such that $\gamma'(0)=m$ and 
$\gamma=h\circ \gamma'$.
\begin{enumerate}
\item   We say that $h$ has the {\it 
path lifting} property if  a lift $\gamma'$ exists for every
$\gamma:[0,1]\to N$ and  $m\in h^{-1}\bigl(\gamma(0)\bigr)$.
  We do not require $\gamma'$ to be unique.  If this holds then every
path $\gamma:[0,1)\to N$ also lifts. 
\item 
We say that $h$  has the {\it arc lifting} (or {\it local
path lifting}) property if the lift $\gamma'$ exists  over some subinterval
$[0,\epsilon]$, where  $\epsilon>0$  may depend on $m\in h^{-1}\bigl(\gamma(0)\bigr)$. 
\item  We say that $h$ has the {\it 
2-point path lifting} property if  given  $\gamma:[0,1]\to N$ and
$m_i\in h^{-1}(\gamma(i))$ for $i=0,1$, there is a lifting $\gamma'$ such that
$\gamma'(0)=m_0$ and $\gamma'(1)=m_1$. 
The  basic example is  a  fiber bundle $h:M\to N$ with path-connected fiber $F$.
\end{enumerate}
Note also that if $N$  is path-connected, $M\neq \emptyset$ and
$h:M\to N$ has the path lifting property then $h$ is surjective.

Lifting constant paths shows that if $h$ has the 2-point path lifting property
then all fibers of $h$ are path-connected.
\end{defn}

The concept of path lifting  occurs most frequently in the topological literature, but
the next result shows that,
from the scheme-theoretic point of view, arc lifting is the most natural property.

\begin{thm} \label{un.op.=.path.lift.loc.thm}
Let  $g:X\to Y$ be a  morphism of $\c$-schemes of finite type. The following are equivalent.
\begin{enumerate}
\item  $g:X\to Y$ is universally open (Definition~\ref{op.uop.defn}).
\item $g(\c):X(\c)\to Y(\c)$ is  open  in the Euclidean topology.
\item $g(\c):X(\c)\to Y(\c)$  has the arc lifting property (\ref{weak.l.p.say}.2).
\end{enumerate}
\end{thm}

Proof.  The equivalence of (\ref{un.op.=.path.lift.loc.thm}.1) and (\ref{un.op.=.path.lift.loc.thm}.2) is established in
Lemma~\ref{vol.op.crit.7}.  

Assume next that $g$ is not  universally open at $x$.
By Theorem~\ref{uop.et.test.thm}, after a suitable \'etale base change, it is not open at $x$. An \'etale base change is a local homeomorphism in the Euclidean topology, thus it does not alter the liftability of arcs.

Thus, by (\ref{vol.op.crit}.2)  there is a curve
$g(x)\in C\subset Y$ such that $f(U)\cap C=\{g(x)\}$. 
Choosing a path in this curve $C$ shows that (\ref{un.op.=.path.lift.loc.thm}.3) $\Rightarrow$ (\ref{un.op.=.path.lift.loc.thm}.2). 
Thus it remains to prove that (\ref{un.op.=.path.lift.loc.thm}.2) $\Rightarrow$ (\ref{un.op.=.path.lift.loc.thm}.3). 
This  done in 3 steps.

If  $g$ is finite, arc lifting is proved in Lemmas~\ref{p.uo.pl.prop.l1}--\ref{p.uo.pl.prop.l2}.
If $X_{g(x)}$ has pure dimension $n$ at $x$, then we take a
relative complete intersection $x\in Z\subset X$ of codimension $n$. 
Then $g|_Z:Z\to Y$ is quasi-finite at $x$  and 
open by (\ref{pd.crit.say}.3). 
An  \'etale base change as in Proposition~\ref{qf.to.f.by.et.thm}
reduces this to the already discussed finite case.
Thus an arc $\gamma:[0,\epsilon]\to Y$
has a lift $\gamma':[0,\epsilon]\to Z$. This is also a lift to $X$.

Finally, 
in general set $n=\dim_x X_{g(x)}$. By (\ref{pd.crit.say}.2)
the restriction of $g$ to $X^{(n)}$ is universally open. The fibers of
$X^{(n)}\to Y$ have  pure dimension $n$, thus $X^{(n)}\to Y$
has the arc lifting property and so does $X$. \qed

\begin{lem}\label{vol.op.crit.7}  
 $g:X\to Y$  be a morphism between $\c$-schemes of finite type.
Then $g$ is universally open at $x$ in the Zariski topology iff 
 $g(\c):X(\c)\to Y(\c)$ is  open  at $x$ in the Euclidean topology.
 \end{lem}

Proof. Assume that $g(\c)$ is open  at $x$ in the Euclidean topology
 and let $x\in U\subset X$ be a Zariski open neighborhood.
Then $U$ is also Euclidean open, hence $g(U)$ contains a 
Euclidean open neighborhood of $g(x)$. Since $g(U)$ is also constructible,
it also contains a Zariski open neighborhood of $x$. Thus
$g$ is  open at $x$ in the Zariski topology.
An \'etale base change is a local homeomorphism in the
 Euclidean topology, hence 
$g$ is also universally open at $x$ in the Zariski topology by Theorem~\ref{uop.et.test.thm}. 

Conversely, assume that $g$ is not  universally open at $x$.
By Theorem~\ref{uop.et.test.thm}, after a suitable \'etale base change, it is not open at $x$.
Thus, by (\ref{vol.op.crit}.1)   this means that there is a curve
$g(x)\in C\subset Y$ such that $f(U)\cap C=\{g(x)\}$.  This 
 shows that $g(\c)$ is not open at $x$ in the Euclidean topology either.
\qed

\begin{lem}  \label{p.uo.pl.prop.l1}
  Let  $g:X\to Y$ be a   finite  morphisms of $\c$-schemes of finite type. Then one can choose triangulations on them such that
  $g(\c):X(\c)\to Y(\c)$ is simplicial.
\end{lem}

Proof.   
Set $n:=\dim Y$, $X_n:=X$ and  $Y_n:=Y$. If $g_r:X_r\to Y_r$ is already defined then let
$U_r\subset Y_r$ be the largest smooth, open subset of pure dimension $r$ such that
$g_r$ is \'etale. Set $Y_{r-1}:=Y_r\setminus U_r$ and $X_{r-1}:= \red g^{-1}(Y_{r-1})$.

Triangulate $Y$ by starting with $Y_0$, extending it to a triangulation of $Y_1$ (possibly adding new vertices), then extending it to a triangulation of $Y_2$ (possibly after a refinement) and so on. At the end we get a triangulation of $Y$ such that the  interior 
$\Delta^\circ$ of every simplex lies in some $ U_{r}$. Thus the preimages of the simplices give a triangulation of $X$ and $g$ is simplicial. \qed

\begin{lem} \label{p.uo.pl.prop.l2}
Let  $M,N$ be connected simplicial complexes and $g:M\to N$ a  proper, open, simplicial map with finite fibers. Then $g$ has the path lifting property.
\end{lem}

Proof. 
We are given a path $\gamma: [0,1]\to N$ and a lifting $\gamma'(0)$. 
Choose a maximal lifting  $\gamma'_c: [0,c)\to N$. (At the begining we allow
$c=0$.)  
First we define $\gamma'(c)$.
 Set $n=\gamma(c)$ and let $m_i\in M$ be the preimages of $n$. Choose an open neighborhood  $n\in U\subset N$ such that $g^{-1}(U)=\cup_i V_i$ is a disjoint union of connected neighborhoods of the $m_i$. Then $\gamma'(c-\epsilon, c)$ is contained in one of them, say $V_j$, and setting $\gamma'(c)=m_j$ is the unique 
continuous extension of  $\gamma':[0,c)\to M$ to $\gamma':[0,c]\to M$. 

In order to extend beyond $c$, we may as well assume that $n:=\gamma(c)$ is a vertex.
Then $m:=\gamma'(c)$ is also a vertex. We can choose the neighborhood $n\in U\subset N$ to be a cone over the link $L_n$.    Then
$m\in V\subset M$ is also a cone over its link $L_m$ and
$g|_{L_m}:L_m\to L_n$ is finite, open and simplicial. 
If $\gamma$ maps  $[c,c+\epsilon]$ to $V$ then we define
$\gamma'$ on $[c,c+\epsilon]$ as follows. If $\gamma(x)=n$ then set
$\gamma'(x)=m$. The complement of this set, denoted by  $B\subset [c,c+\epsilon]$,
is a countable union of open intervals and $\gamma$ maps $B$ to
$$
U\setminus\{n\}\sim (0,1)\times L_n.
$$
Correspondingly we can write  $\gamma|_B=(\alpha_B, \gamma_B)$ where
$\gamma_B$ maps $B$ to $L_n$.
Since $g|_{L_m}:L_m\to L_n$ is simplicial, by induction on the dimension,
$\gamma_B$  lifts to $\gamma'_B:B\to L_m$. (As we noted in (\ref{weak.l.p.say}.1), open paths also lift.)
Then set
$\gamma'|_B:=(\alpha_B, \gamma'_B)$.\qed

\medskip
Next we consider the scheme-theretic description of the 
2-point path lifting property.

\begin{thm} \label{un.op.=.path.lift.thm}
Let  $g:X\to Y$ be a   morphism of $\c$-schemes of finite type. The following are equivalent.
\begin{enumerate}
\item  $g:X\to Y$ is universally open,  surjective and has connected fibers.
\item $g(\c):X(\c)\to Y(\c)$ is  open in the Euclidean topology,  surjective  and has connected fibers.
\item $g(\c):X(\c)\to Y(\c)$ is    surjective,   has connected fibers and satisfies the arc lifting property (\ref{weak.l.p.say}.2).
\item $g(\c):X(\c)\to Y(\c)$  has the 2-point path lifting property (\ref{weak.l.p.say}.3). 
\end{enumerate}
\end{thm}

Proof. The equivalence of (1), (2) and (3) follows from 
Theorem~\ref{un.op.=.path.lift.loc.thm}, once we note that 
 a $\c$-scheme of finite type is connected in the Zariski topology iff it is 
connected in the Euclidean topology. (This is usually proved using Chow's theorem, but one can use Bertini's hyperplane section theorem to reduce it to the 1-dimensional case, which was known to Riemann.)

We noted  at the end of Definition~\ref{weak.l.p.say} that if (4) holds then $g$ is surjective and has connected fibers.
Thus (4) implies (1--3) and it remains to show that (1--3) imply (4).

We have $\gamma:[0,1]\to Y(\c)$ and liftings $m_0$ of $\gamma(0)$
and $m_1$ of $\gamma(1)$.
Assume that we already have a lifting $\gamma'$ defined on $[0,c]$.
By (3) we have a lifting  $\gamma''_c:[c, c+\epsilon]\to X(\c)$ such that
$\gamma''_c(c)=\gamma'(c) $. The concatenation of $\gamma'$ on $[0,c]$ with $\gamma''_c$
 defines $\gamma'$ on $[0,c+\epsilon]$.

 Assume next that we have  $\gamma'$ defined on $[0,c)$.
   We do not claim that it extends to  $[0,c]$, but we show that
   the restriction of $\gamma'$ to $[0,c-\epsilon]$ extends to  $\gamma'':[0,c+\epsilon]$
   for some $\epsilon>0$. Moreover, if $c=1$ then we can also arrange that
$\gamma''(1)=m_1$. For other values of $c$ pick any lifting
$m_c$ of $\gamma(c)$. Applying (3) gives a lifting
 \begin{enumerate}\setcounter{enumi}{4}
     \item $\gamma''_c:[c-\epsilon, c+\epsilon]\to X$ such that 
 $\gamma''_c(c)=m_c$. 
 \end{enumerate}
Starting from  $X\to Y$, we get a stratification  $\{Y_i\}$ of $Y$  as in
(\ref{top.triv.start}.1). 
   Applying Lemma~\ref{p.uo.pl.prop.l2} to  $\amalg_iY_i\to Y$  
   we get a triangulation of  $Y$ with $\gamma(c)$ as a vertex.
Let $L$ be the smallest subcomplex that contains
$\gamma\bigl((c-\delta, c)\bigr)$ for some $\delta>0$ and let 
 $\Delta\subset L$ be a maximal dimensional simplex with interior $\Delta^\circ$.
Then $\gamma^{-1}(\Delta^\circ)$ is open in $(c-\delta,c)$ and its closure contains $c$. 
In particular, we see that 
   \begin{enumerate}
     \item[(a)] there is a homeomorphism  $\beta:g^{-1}(\Delta^\circ)\sim \Delta^\circ\times M$  that   commutes with projection to  $\Delta^\circ$ for some $M$ and
     \item[(b)] there  are  $0<\eta_1<\eta_2<\epsilon$ such that
$\gamma$ maps  $[c-\eta_2, c-\eta_1]$ to $\Delta^\circ$.
   \end{enumerate}
As we noted in (\ref{weak.l.p.say}.3),
   the restriction of $\gamma$ to $[c-\eta_2, c-\eta_1]$
   has a lifting $\bar\gamma$ such that
   $$
   \bar\gamma(c-\eta_2)=\gamma'(c-\eta_2)\qtq{and}
   \bar\gamma(c-\eta_1)=\gamma''_c(c-\eta_1).
   $$
   Thus the concatenation of  $\gamma'$ on $[0,c-\eta_2]$,
   $\bar\gamma$ on  $[c-\eta_2,  c-\eta_1]$ and
   $\gamma''_c$ on  $[c-\eta_1,  c+\epsilon]$
   defines a lifting  of $\gamma$ over $[0,c+\epsilon]$. \qed

\begin{say}[Stratification of maps]\label{top.triv.start}
Let  $g:X\to Y$ be a morphism of varieties over $\c$.
As in \cite[Sec.I.1.7]{gm-book},
$Y$ has a stratification by closed subvarieties
$$
Y=Y_n\supset Y_{n-1}\supset \cdots \supset  Y_0
\eqno{(\ref{top.triv.start}.1)}
$$
where $Y_i\setminus Y_{i-1}$ has pure dimension $i$ and 
each  $$
g^{-1}\bigl(Y_i\setminus Y_{i-1}\bigr)\to \bigl(Y_i\setminus Y_{i-1}\bigr)
$$
is a topologically locally trivial fiber bundle.
(There may be different fibers over different connected components.)
Thus path lifting could fail only at points where a path moves from one stratum to another. 
\end{say}

\medskip

Let  $g:X\to Y$ be a quasi-finite, universally open morphism. 
Using Lemma~\ref{p.uo.pl.prop.l2} we see that  $g(\c):X(\c)\to Y(\c)$  has the  path lifting property 
iff $g$ is proper. 
However, as shown by Example~\ref{prop.pl.exmps}, 
the  path lifting property does not seem to have an equivalent scheme-theoretic version for morphisms with positive dimensional fibers. Nonetheless,  the following sufficient  condition is quite natural and useful. 

\begin{thm} \label{un.op.=.path.lift.thm.2}
Let  $g:X\to Y$ be a  proper, universally open, pure dimensional  morphism of $\c$-schemes of finite type. Then  $g(\c):X(\c)\to Y(\c)$  has the  path lifting property 
\end{thm}

Proof.  As  (\ref{basic.props.say}.8) shows, we can not use Stein factorization to reduce Theorem~\ref{un.op.=.path.lift.thm.2}
directly to Theorem~\ref{un.op.=.path.lift.thm}, but a suitable modification of the proof will work.

We follow the proof of (\ref{un.op.=.path.lift.thm}.3) $\Rightarrow$ (\ref{un.op.=.path.lift.thm}.4),  but we need to make a different choice of
$\gamma''_c$ in  (\ref{un.op.=.path.lift.thm}.5). 
Note that $\gamma'\bigl([c-\eta_1, c-\eta_2]\bigr)$ is contained in a 
connected component  $\Delta^\circ\times M_1\subset \Delta^\circ\times M$, and 
everything works as before if we can ensure that
$\gamma''_c\bigl([c-\eta_1, c-\eta_2]\bigr)$
is also contained in $\Delta^\circ\times M_1$. 

To do this, let $X\to Y_1\to Y$ be the Stein factorization.
The choice of the connected component  $\Delta^\circ\times M_1\subset \Delta^\circ\times M$ defines a section $\sigma: \Delta^\circ\to
 Y_1$. By Lemma~\ref{un.op.=.path.lift.thm.2.lem}, after an \'etale base change
we may assume to have a relative complete intersection $Z\subset X$
such that  the  induced map  $g_1:Z\to Y_1$ is finite  and surjective. 
In particular, there is an $m_z\in Z$ such that $g_1(m_z)=\sigma\bigl(\gamma(c-\eta_2)\bigr)$. Since $Z\to Y$ is finite and open, it satisfies the path lifting property by Lemma~\ref{p.uo.pl.prop.l2}. Thus the restriction of $\gamma$ to $[c-\eta_2, c]$ has a lifting $\gamma''_c$ such that
$\gamma''_c(c-\eta_2)=m_z$. By construction $\gamma''_c(c-\eta_2)$ and
$\gamma'(c-\eta_2)$ are in the same commected component of the fiber and
the rest of the proof now works as before. \qed

\begin{lem} \label{un.op.=.path.lift.thm.2.lem}
Let  $g:X\to Y$ be a  proper, universally open  morphism of  finite type and of relative dimension $n$.   Let $X\to  Y_1\to Y$ denote its Stein factorization.
Let $y\in Y$ be a point and $P\subset X_y$ a finite subset that has nonempty intersection with every irreducible component of $X_y$. 
Let $P\subset Z\subset X$ be a relative complete intersection of codimension $n$. Then there is an \'etale neighborhood  $(y',Y')\to  (y, Y)$ such that
the induced morphism  $Z'\to  Y'_1$ is finite and surjective. 
\end{lem}

Proof. By construction $Z\to Y$ is quasi-finite, hence  there is  
an \'etale neighborhood  $(y',Y')\to  (y, Y)$ such that
 $Z'\to Y'$ is finite.  Thus  $Z'\to  Y'_1$ is also finite. We need to prove that it is  surjective.  

Let $\tilde Y'_1$ be the disjoint union of the irreducible components of $Y'_1$. By base change we get $\tilde X'_1\to \tilde Y'_1$,
$\tilde P'_1\subset \bigl(\tilde X'_1)_{y'}$ and
$\tilde P'_1\subset \tilde Z'_1\subset \tilde X'_1$. 
An irreducible component of $\bigl(\tilde X'_1)_{y'}$ has pure dimension $n$,
thus its image in $X_y$ is an irreducible component. Thus $\tilde P'_1$ has nonempty intersection with every irreducible component of $\bigl(\tilde X'_1)_{y'}$. Thus $\tilde Z'_1$ dominates $\tilde Y'_1$ and hence
$Z'$ dominates $Y'_1$.  Since $Z'\to  Y'_1$ is  finite, this implies that it is surjective. \qed

\medskip

Next we establish the Euclidean version of Corollary~\ref{fin.ind.cor}.

\begin{cor} \label{fin.ind.cor.euc}
Let $g:(x,X)\to (y,Y)$ be a proper, universally open morphism of pointed,  connected $\c$-schemes.  Then
$\im[\pi_1\bigl(X(\c), x\bigr)\to \pi_1\bigl(Y(\c), y\bigr)]$ has finite index in
$\pi_1\bigl(Y(\c), y\bigr)$. 
\end{cor}

Proof. Let $n$ denote the maximal fiber dimension.
Then $g^{(n)}: X^{(n)}\to Y$ is also proper and  universally open by (\ref{pd.crit.say}.2).
We may choose  $x\in X^{(n)}$ and then $g^{(n)}_*$ factors as
$$
\pi_1\bigl(X^{(n)}(\c), x\bigr)\to\pi_1\bigl(X(\c), x\bigr)\to  \pi_1\bigl(Y(\c), y\bigr).
$$
Thus it is enough to show that $\im[\pi_1\bigl(X^{(n)}(\c), x\bigr)\to \pi_1\bigl(Y(\c), y\bigr)]$ has finite index in
$\pi_1\bigl(Y(\c), y\bigr)$. 
The advantage is that   $g(\c):X^{(n)}(\c)\to Y(\c)$ has the
path lifting property by Theorem~\ref{un.op.=.path.lift.thm.2}.
The rest follows from (\ref{fin.ind.cor.euc}.1). 
\medskip

{\it Claim \ref{fin.ind.cor.euc}.1.} Let $h:M\to N$  be a continuous map of path-connected topological spaces  that has the  path lifting property. Assume that  $h^{-1}(n)$ has finitely many path-connected components for some $n\in N$. Then the image of $\pi_1(M,m)\to \pi_1(N,n) $ has finite index in $\pi_1(N,n) $.
\medskip

Proof.   Every loop $\gamma$ starting and ending in $n$ lifts to a path that 
starts at $m$ and ends in $h^{-1}(n)$. If 2 loops  $\gamma_1, \gamma_2$ end
at the same path-connected component then $\gamma_1^{-1}\gamma_2$ lifts to a loop on $M$. This
shows that the index of the image of $\pi_1(M,m)\to \pi_1(N,n) $
is bounded by the number of path-connected components of the fiber. \qed

\subsection*{Examples}{\ }

The first example shows that
Theorem~\ref{un.op.=.path.lift.thm} does not have a pointwise version.

\begin{exmp}\label{op.not.op.exmp.2}
Let $Y:=(x^2=y^2z)$ be the pinch point as in Example~\ref{op.not.op.exmp}.1
with normalization $p:\bar Y\to Y$. We saw that $p$ is universally open 
at the origin. However, 
it does not have the arc lifting property at the origin.

To see this consider the real curves $\gamma^{\pm}:t\mapsto \bigl(t,\pm t^2\sin(t^{-1})\bigr)\subset \r^2$. Note that  $p\circ \gamma^+$ and
$p\circ \gamma^-$ intersect at the points  $t=(m\pi)^{-1}$. Thus
 the arc
$$
\gamma(t):=
\left\{
\begin{array}{l}
p\circ \gamma^+(t)\qtq{if}  \gamma^+(t)\geq 0\qtq{and}\\
p\circ \gamma^-(-t)\qtq{if}  \gamma^-(-t)\geq 0
\end{array}
\right.
$$
has no  lifting.
\end{exmp}

\begin{exmp}[Path lifting and properness]\label{prop.pl.exmps} 
It is natural to hope that   arc lifting plus properness should imply path lifting, but this is not the case.
\medskip

(\ref{prop.pl.exmps}.1)   Let $X$ be obtained from
$X_1:=B_{(0,0)}\c^2$ and of $X_2:=\p^1\times \c^2$ by identifying
the exceptional divisor of the blow-up with  $\p^1\times \{(0,0)\}$. 
The projection $X\to \c^2$ is universally open, and so is
$X_2\to \c^2$, but $X_1\to \c^2$ is not even open. 
The path 
$$
t\mapsto  \bigl(1-t, (1-t)\sin\bigl(\tfrac1{1-t}\bigr)\bigr)
$$
does not lift to $B_{(0,0)}\c^2$. Thus 
 $X\to \c^2$ is universally open and proper, it has the
arc lifting property but not the path lifting property.
\medskip

In the next 2 variants of the above construction, gluing of 2 irreducible components has the opposite effect on path lifting.
\medskip

(\ref{prop.pl.exmps}.2) 
Set $Y_1:=(ux=vy)\subset \p^2_{xyz}\times \c^2_{uv}$.
Note that $\Sigma:=(x=y=0)\subset Y_1$ is a section of
$\pi_2:Y_1\to \c^2_{uv}$. Thus every path starting in
$\c^2_{uv}\setminus\{(0,0)\}$ can be lifted to $Y_1$ but
$Y_1\to \c^2_{uv}$ does  not have the path lifting property.

Let $Y$ be obtained from
$Y_1$ and of $Y_2:=\p^2\times \c^2_{uv}$ by identifying
the fibers over the origin.   Then $p:Y\to \c^2_{uv}$ has the path lifting property.
\medskip

(\ref{prop.pl.exmps}.3) 
Let $Z_1$ be obtained from  $\p^1\times \c^2_{uv}$ by blowing up
 $\p^1\times \{(0,0)\}$. The projection to  $\c^2_{uv}$ factors as $Z_1\to B_{(0,0)}\c^2\to \c^2$, thus not all paths starting in $\c^2_{uv}\setminus\{(0,0)\}$ can be lifted to $Z_1$.
Let $Z$ be obtained from
$Z_1$ and of $Z_2:=\p^1\times \p^1\times \c^2_{uv}$ by identifying
the fibers over the origin.  Thus 
 $Z\to \c^2$ is universally open and proper, it has the
arc lifting property but not the path lifting property.
\end{exmp}

\begin{ack} 
Partial  financial support    was provided  by  the NSF under grant numbers
 DMS-1362960 and  DMS-1440140 while the author was in residence at
MSRI during the Spring 2019 semester.
\end{ack}

%\bibliography{refs}
\def\cprime{$'$} \def\cprime{$'$} \def\cprime{$'$} \def\cprime{$'$}
  \def\cprime{$'$} \def\dbar{\leavevmode\hbox to 0pt{\hskip.2ex
  \accent"16\hss}d} \def\cprime{$'$} \def\cprime{$'$}
  \def\polhk#1{\setbox0=\hbox{#1}{\ooalign{\hidewidth
  \lower1.5ex\hbox{`}\hidewidth\crcr\unhbox0}}} \def\cprime{$'$}
  \def\cprime{$'$} \def\cprime{$'$} \def\cprime{$'$}
  \def\polhk#1{\setbox0=\hbox{#1}{\ooalign{\hidewidth
  \lower1.5ex\hbox{`}\hidewidth\crcr\unhbox0}}} \def\cdprime{$''$}
  \def\cprime{$'$} \def\cprime{$'$} \def\cprime{$'$} \def\cprime{$'$}
\providecommand{\bysame}{\leavevmode\hbox to3em{\hrulefill}\thinspace}
\providecommand{\MR}{\relax\ifhmode\unskip\space\fi MR }
\providecommand{\MRhref}[2]{%
  \href{http://www.ams.org/mathscinet-getitem?mr=#1}{#2}
}
\providecommand{\href}[2]{#2}

\bigskip

\noindent  Princeton University, Princeton NJ 08544-1000

\email{kollar@math.princeton.edu}

\end{document}